\RequirePackage{fix-cm}
\documentclass[smallextended]{svjour3}       
\smartqed  
\usepackage{graphicx}
\journalname{MATRI}
\begin{document}

\title{Bj\"orck-Pereyra-type methods and total positivity
}


\author{Jos\'e-Javier Mart{\'\i}nez
}


\institute{J.-J. Mart{\'\i}nez \at
              Departamento de F{\'\i}sica y Matem\'aticas, \\
              Universidad de Alcal\'a, \\
              Alcal\'a de Henares, Madrid 28871, Spain\\
              \email{jjavier.martinez@uah.es}           
}

\date{Received: date / Accepted: date}

\maketitle

\begin{abstract}
The approach to solving linear systems with structured matrices by means of the bidiagonal factorization of the inverse of the coefficient matrix is first considered, the starting point being the classical Bj\"orck-Pereyra algorithms for Vandermonde systems, published in $1970$ and carefully analyzed by Higham in $1987$. The work of Higham showed the crucial role of total positivity for obtaining accurate results, which led to the generalization of this approach to totally positive Cauchy, Cauchy-Vandermonde and generalized Vandermonde matrices.

Then, the solution of other linear algebra problems (eigenvalue and singular value computation,  least squares problems) is addressed, a fundamental tool being the bidiagonal decomposition of the corresponding matrices. This bidiagonal decomposition is related to the theory of Neville elimination, although for achieving high relative accuracy the algorithm of Neville elimination is not used. Numerical experiments showing the good behaviour of these algorithms when compared with algorithms which ignore the matrix structure are also included.

\end{abstract}

\section{Introduction}
\label{intro}

The second edition of the Handbook of Linear Algebra \cite{Handbook}, edited by Leslie Hogben, is substantially expanded from the first edition of $2007$ and, in connection with our work, it contains a new chapter by Michael Stewart entitled {\it Fast Algorithms for Structured Matrix Computations} (chapter $62$) \cite{Stewart}. This chapter includes, as Section $62.7$, the subject of fast algorithms for Vandermonde systems. Among these algorithms the author incorporates the Bj\"orck-Pereyra algorithm for solving Vandermonde linear systems, indicating the relationship with polynomial interpolation by using the Newton basis, and the interpretation of this process as a factorization in terms of bidiagonal matrices.

Section $62.7$ also recalls the high relative accuracy of the Bj\"orck-Pereyra algorithm when the (nonnegative) nodes are ordered in increasing order, a fact observed in the error analysis presented by N. Higham in \cite{Higham}.

In connection with {\it high relative accuracy}, the book \cite{Handbook} also contains a chapter (no. $59$), due to Z. Drma{\v c}, entitled {\it Computing Eigenvalues and Singular Values to High Relative Accuracy} (already available in the first edition as chapter $46$) \cite{Drmac}. In this chapter the author includes the references to the work of Demmel and Kahan \cite{DKah} and of Fernando and Parlett \cite{FP} on the computation of singular values of bidiagonal matrices, fundamental references for the subsequent work of Demmel and Koev (see \cite{DK}, \cite{KOEV05}, \cite{KOEV07}).

\medskip

We observe the same situation in two other recent relevant books on numerical linear algebra: the fourth edition of Golub and Van Loan's {\it Matrix Computations} \cite{GolV} and the new book by A. Bj\"orck {\it Numerical Methods in Matrix Computations} \cite{Bjorck}. In both works the authors pay attention to the description of the Bj\"orck-Pereyra algorithm for solving Vandermonde linear systems, in the corresponding chapters devoted to linear system solving: chapter $4$ ({\it Special Linear Systems} in \cite{GolV}, and chapter $1$ ({\it Direct Methods for Linear Systems}) in \cite{Bjorck}. Other separate chapters from both books are devoted to eigenvalue and singular value problems.

Consequently, to the reader these problems appear as very different problems: linear system solving, or eigenvalue and singular value problems. The purpose of the present work is to show, in the important cases where total positivity is present, the unifying role of the bidiagonal decompositions of the corresponding matrices to develop a first stage in algorithms for solving various linear algebra problems with high relative accuracy.

The analysis of the Bj\"orck-Pereyra algorithms presented by Higham in \cite{Higham} showed (among many other things) the important role of total positivity in the accuracy of the algorithms. So this property will be of fundamental importance in our account.

Let us recall that, according to the classical respectable definition, a matrix is said to be {\sl totally positive} if all its minors are nonnegative, and when all the minors are positive the matrix is called {\sl strictly totally positive} \cite{MM18}. More recently, matrices with all their minors nonnegative have been called {\sl totally nonnegative} matrices, this being mathematically more precise.

Recent books on this subject are \cite{FJ,Pinkus}. These monographs cover many aspects of the theory and applications of totally positive matrices and, although they do not develop the topic of accurate computations with this type of matrices, they contain useful references to the work of Demmel and Koev in this field.

The classical Bj\"orck-Pereyra algorithms date back to $1970$, and now we want to highlight the important work of Demmel and Koev in the development of these algorithms with high relative accuracy. In \cite{DK} Demmel and Koev called some of these methods Bj\"orck-Pereyra-type methods, acknowledging in this way the importance of the work of Bj\"orck and Pereyra, and in the present work we want to carry out a necessary explanation of the analogies and differences between Bj\"orck-Pereyra-type methods and related methods, the main analogy being the role of bidiagonal factorization of matrices.

One key fact is that while Bj\"orck-Pereyra algorithms for linear systems were related to the factorization of the inverse $A^{-1}$ of the coefficient matrix, the use of the bidiagonal factorization of $A$ considered in the work of Koev allowed to extend the accurate computations to various linear algebra problems different from linear system solving.

The rest of the paper is organized as follows. In Section $2$ the problem of linear system solving is considered, a problem to which Bj\"orck-Pereyra algorithms and their generalizations were devoted. Section $3$ reviews Neville elimination and the bidiagonal factorization associated to it, a key theoretical tool for the development of algorithms. The extension of these algorithms based on a bidiagonal factorization to the problems of eigenvalue and singular value computation is addressed in Section $4$, while Section $5$ considers the extension to the rectangular case, which allows solving least squares problems. The brief Section $6$ considers the class of totally positive tridiagonal matrices, and finally Section $8$ is devoted to conclusions.

\section{Linear system solving: the Vandermonde case and some extensions}

As recently recalled in \cite{Bjorck}, it was in $1970$ when Bj\"orck and Pereyra showed that the Newton-Horner algorithm for solving a Vandermonde linear system related to polynomial interpolation could be expressed in terms of a factorization of the inverse of the Vandermonde matrix as a product of diagonal and lower bidiagonal matrices. In their paper \cite{BP} those authors included, after the numerical experiments, the following sentence:

\medskip
\textit{It seems as if at least some problems connected with Vandermonde systems, which traditionally have been considered too ill-conditioned to be attacked, actually can be solved with good precision}.
\medskip

Many years later, at the beginning of his brilliant contributions in the field on numerical  linear algebra, Higham gave in \cite{Higham} an analysis of Bj\"orck-Pereyra algorithms and indicated that when the interpolation nodes are nonnegative and ordered in increasing order then the corresponding Vandermonde matrix is \textit{totally positive}. If, in addition, the components of $f$ alternate in sign, then there is no cancellation and high relative accuracy is obtained.

\medskip

We follow the clear presentation of \cite{GolV} and \cite{Bjorck} to show the matrix interpretation of the Bj\"orck-Pereyra algorithm to solve the dual Vandermonde system $V^Ta = f$, i.e. the linear system corresponding to polynomial interpolation. These authors call {\it Vandermonde matrix} to the transpose $V$ of that matrix, whose first row is $(1,1,\ldots,1)$.

As we read in Section $4.6$ of \cite{GolV} and in subsection $1.8.3$ of \cite{Bjorck}, the solution $a$ of $V^Ta = f$ is given as

$$
a = V^{-T} f = L^T U^T f =
$$

$$
= L_0(x_0)^T \cdots L_{n-1}(x_{n-1})^T D_{n-1}^{-1} L_{n-1}(1) \cdots D_{0}^{-1} L_{0}(1),
$$
where $x_0,\ldots,x_n$ are the interpolation nodes, and $V = V(x_0,\ldots,x_n)$ is the corresponding Vandermonde matrix.

For $i=0,\ldots,n-1$, $D_i$ are diagonal matrices, while $L_{i}(x_i)$ and $L_i(1)$ are lower triangular bidiagonal matrices with ones on the main diagonal. In section $3$ we will show with an example the precise structure of those matrices and we will show the differences between this factorization and the factorization related to Neville elimination.

Let us observe that we can identify two stages in the algorithm, since $a = L^T U^T f = L^T c$, where $c = U^T f$. This vector $c$ is the vector of divided differences, i.e. the coefficients of the interpolating polynomial in the {\it Newton basis}. The second stage ($a = L^Tc$) corresponds to the change of basis from the Newton basis to the monomial basis.

\medskip

This approach was later applied by Boros, Kailath and Olshevsky to design a Bj\"orck-Pereyra-type algorithm for solving Cauchy linear systems. As we read in \cite{BKO}, ``\textit{for the class of totally positive Cauchy matrices the new algorithm is forward and backward stable, producing a remarkable high relative accuracy. In particular, Hilbert linear systems, often considered to be too ill-conditioned to be attacked, can be rapidly solved with high precision}''.

Other generalization (in this case to confluent Vandermonde-like systems) was presented by Higham (see \cite{HIGSIAM} and Chapter $22$ of his book \cite{HIG}), but in general the bidiagonal structure of the matrices involved in the factorization is lost. More recently a new generalization of Bj\"orck-Pereyra algorithms (now to the case of Szeg\"o-Vandermonde matrices) was introduced in \cite{Bella}, but the authors admit that some matrices involved in the factorization are not sparse (i.e. they are far from being bidiagonal).

A bidiagonal factorization of the inverse of the coefficient matrix (in this case for a Cauchy-Vandermonde matrix) related to Neville elimination was presented in \cite{MP98}, a pioneering paper on this subject which has recently been completed in the light of the work of Koev (see \cite{MMP}).
Also, in \cite{DK} the authors extend the Bj\"orck-Pereyra algorithm to solve totally positive generalized Vandermonde linear systems $G y = b$ by computing a bidiagonal decomposition of $G^{-1}$ and computing $y = G^{-1}b$. The authors indicate that this decomposition, related to Neville elimination \cite{GP96}, reveals the total positivity of $G$ and yields a Bj\"orck-Pereyra-type algorithm for the solution of $G y = b$.

They also comment that when $b$ has an alternating sign pattern then their algorithm is subtraction-free, so that the solution $y$ is very accurate, the same fact observed by Higham in \cite{Higham} when analyzing the Bj\"orck-Pereyra algorithm. This is a consequence of the fact that the inverse of a totally positive matrix has a checkerboard sign pattern.

In their work \cite{DDHK} Demmel and co-workers have briefly indicated the importance of these Bj\"orck-Pereyra-type methods for solving linear systems with several classes of totally positive matrices (Vandermonde, Cauchy, Pascal, Cauchy-Vandermonde, generalized Vandermonde, Bernstein-Vandermonde,...). The main fact, however, is that when the bidiagonal factorization of the coefficient matrix $A$ is related to \textit{Neville elimination}, then there are exactly $n^2$ independent nonnegative parameters in the bidiagonal decomposition which are stored in a matrix $B = \mathcal{BD(A)}$ (bidiagonal decomposition of $A$). Then, starting from an accurate $\mathcal{BD(A)}$, virtually all linear algebra with totally nonnegative matrices can be performed accurately by using algorithms due to P. Koev \cite{KOEV}.

The miracle of the work of Koev was to show that, while all these Bj\"orck-Pereyra-type methods pay attention to the factorization of the inverse $A^{-1}$ of the coefficient matrix $A$ (which seems natural when solving a linear system $A x = f$), if one uses the (unique) bidiagonal decomposition of $A$ related to Neville elimination various linear algebra problems can additionally be solved.

Therefore, our following section is devoted to present this bidiagonal decomposition of $A$, and to show the analogies and differences between classical Bj\"orck-Pereyra algorithms and algorithms related to Neville elimination when dealing with totally positive matrices.

\section{Bidiagonal factorization and Neville elimination}

As  recently recalled by J. M. Pe\~na in \cite{PenaS}, the work by Gasca and M\"uhlbach on the connection between interpolation formulas and elimination techniques made clear that what they named {\it Neville elimination} had special interest for totally positive matrices (see, for instance, \cite{MG87}). Later on, in \cite{GP92,GP94,GP96} Gasca and Pe\~na greatly improved the previous results on Neville elimination an total positivity, and their work was one of the starting points for the work of Demmel and Koev in developing algorithms which, starting from the appropriate bidiagonal decomposition associated to Neville elimination, made possible the  accurate solution of many linear algebra problems for totally positive matrices, including eigenvalue and singular value computation (see \cite{DK,KOEV05,KOEV07}).

Neville elimination characterizes nonsingular totally positive matrices, according to the results of \cite{GP96} recalled as Theorem $8$ in \cite{PenaS}:

\medskip
\begin{theorem} A square matrix $A$ is nonsingular totally positive if and only if the Neville elimination of $A$ and $A^T$ can be performed without row exchanges, all the multipliers of the Neville elimination of $A$
and $A^T$ are nonnegative and all the diagonal pivots of the Neville
elimination of $A$ are positive.
\end{theorem}

\medskip

The bidiagonal decomposition of $A$ and of its inverse $A^{-1}$ are given in the following theorems (see \cite{GP92,GP94,GP96,KOEV05,KOEV07,PenaS}):

\begin{theorem} Let $A$ be a nonsigular totally positive matrix of size $n \times n$. Then $A^{-1}$ admits a factorization in
the form
$$A^{-1}=G_1 \cdots G_{n-1}D^{-1}F_{n-1}\cdots F_1, $$
where
$G_i$ ($1\le i\le n-1$) are $n \times n$ upper triangular bidiagonal matrices, $F_i$ ($i=1,\ldots,n-1$) are $n \times n$ lower triangular bidiagonal matrices, and $D$ is a diagonal matrix of order $n$.

\end{theorem}

\medskip

\begin{theorem} Let $A$ be a nonsingular totally positive matrix of size $n \times n$. Then $A$ admits a factorization in
the form
$$A= {\bar F}_{n-1} \cdots {\bar F}_1 D {\bar G}_1 \cdots {\bar G}_{n-1}$$
where ${\bar F}_i$ ($1\le i\le n-1$) are $n \times n$ lower triangular bidiagonal matrices,
${\bar G}_i$ ($1\le i\le n-1$) are $n \times n$ upper triangular bidiagonal matrices, and $D$ is a diagonal matrix of order $n$.

\end{theorem}

\medskip

The relationship between the process of Neville elimination and the off-diagonal entries of the bidiagonal matrices (all their diagonal entries are equal to one) and the diagonal entries of $D$ can be seen in  \cite{PenaS} , and in \cite{MM13} for the particular case of Bernstein-Vandermonde matrices.

\medskip

\begin{remark}
\rm The algorithm {\tt TNBD} in \cite{KOEV} computes the matrix denoted
as $\mathcal{BD}(A)$ in \cite{KOEV07}, which represents the {\it bidiagonal
decomposition} of $A$. But it is a remarkable fact that the same
matrix $\mathcal{BD}(A)$ also serves to represent the bidiagonal
decomposition of $A^{-1}$. The algorithm computes $\mathcal{BD}(A)$
by performing {\it Neville elimination} on $A$, which involves true
subtractions, and therefore does not guarantee high relative
accuracy.
\end{remark}

\medskip

As we have seen, the matrix $\mathcal{BD}(A)$ serves to have a new parameterization of a totally nonnegative matrix $A$. Although it does not guarantee high relative accuracy, the algorithm {\tt TNBD} of the package {\tt TNTool} \cite{KOEV} computes $\mathcal{BD}(A)$ starting from $A$. But, remarkably, we also have the other way: starting from $ B = \mathcal{BD}(A)$, the algorithm {\tt TNExpand} of \cite{KOEV} computes (now with high relative accuracy) the matrix $A$. In addition, the algorithm {\tt TNInverseExpand} of Marco and Mart{\'\i}nez (also included in \cite{KOEV}) computes, starting from $ B = \mathcal{BD}(A)$, the inverse matrix $A^{-1}$.

The following example illustrates this facts, and the accuracy obtained by using the algorithm {\tt TNInverseExpand}.

{\it Example.} Starting from

$$
B = \left(
      \begin{array}{cccc}
        16 & 3 & 12 & 13 \\
        & & & \\
        5 & 10 & 11 & 8 \\
        & & & \\
        9 & 6 & 7 & 12 \\
        & & & \\
        4 & 15 & 14 & 1 \\
      \end{array}
    \right),
$$
(the magic square appearing in D\"urer's engraving {\it Melencolia I}) we obtain, by means of the sentence {\tt A = TNExpand(B)}, the following matrix (which we can call the D\"urer totally positive matrix):

$$
A = \left(
      \begin{array}{cccc}
        16 & 48 & 96 & 1248 \\
        & & & \\
        80 & 250 & 610 & 8810 \\
        & & & \\
        720 & 2310 & 6277 & 94941 \\
        & & & \\
        2880 & 10140 & 37011 & 617764 \\
      \end{array}
    \right).
$$

This means that the bidiagonal decomposition of $A$ is

$$
A = \left(
      \begin{array}{cccc}
        1 & 0 & 0 & 0 \\
        & & & \\
        0 & 1 & 0 & 0 \\
        & & & \\
        0 & 0 & 1 & 0 \\
        & & & \\
        0 & 0 & 4 & 1 \\
      \end{array}
    \right)
    \left(
      \begin{array}{cccc}
        1 & 0 & 0 & 0 \\
        & & & \\
        0 & 1 & 0 & 0 \\
        & & & \\
        0 & 9 & 1 & 0 \\
        & & & \\
        0 & 0 & 15 & 1 \\
      \end{array}
    \right)
    \left(
      \begin{array}{cccc}
        1 & 0 & 0 & 0 \\
        & & & \\
        5 & 1 & 0 & 0 \\
        & & & \\
        0 & 6 & 1 & 0 \\
        & & & \\
        0 & 0 & 14 & 1 \\
      \end{array}
    \right)
    \left(
      \begin{array}{cccc}
        16 & 0 & 0 & 0 \\
        & & & \\
        0 & 10 & 0 & 0 \\
        & & & \\
        0 & 0 & 7 & 0 \\
        & & & \\
        0 & 0 & 0 & 1 \\
      \end{array}
    \right)
$$

$$
\left(
      \begin{array}{cccc}
        1 & 3 & 0 & 0 \\
        & & & \\
        0 & 1 & 11 & 0 \\
        & & & \\
        0 & 0 & 1 & 12 \\
        & & & \\
        0 & 0 & 0 & 1 \\
      \end{array}
    \right)
    \left(
      \begin{array}{cccc}
        1 & 0 & 0 & 0 \\
        & & & \\
        0 & 1 & 2 & 0 \\
        & & & \\
        0 & 0 & 1 & 8 \\
        & & & \\
        0 & 0 & 0 & 1 \\
      \end{array}
    \right)
    \left(
      \begin{array}{cccc}
        1 & 0 & 0 & 0 \\
        & & & \\
        0 & 1 & 0 & 0 \\
        & & & \\
        0 & 0 & 1 & 13 \\
        & & & \\
        0 & 0 & 0 & 1 \\
      \end{array}
    \right)
$$

Although of size $4 \times 4$, the condition number of this matrix is $\kappa_2 (A) = 1.4 e+11$, which means the computation of the inverse (which we call $AIm$) by using MATLAB will suffer the effect of the high condition number. In fact, by using the exact inverse (which we call $AIe$) computed by Maple we find that the relative error computed in the spectral norm is

$$
\frac{\| AIm - AIe \|_2}{\| AIe \|_2} = 2.8 e-10 ,
$$
while for the inverse matrix $AIk$ computed by using the algorithm {\tt TNInverseExpand(B)} we have full accuracy:

$$
\frac{\| AIk - AIe \|_2}{\| AIe \|_2} = 1.2 e-16.
$$

\bigskip

Now let us show the differences (also the analogies) between the bidiagonal factorization associated with Neville elimination (stored in $B = \mathcal{BD}(A)$) and the factorization corresponding to the Bj\"orck-Pereyra algorithm, by using as example the Vandermonde matrix of order $n=4$ with nodes $2, 3, 5, 8$:

$$
A = \left(
      \begin{array}{cccc}
        1 & 2 & 4 & 8 \\
        & & & \\
        1 & 3 & 9 & 27 \\
        & & & \\
        1 & 5 & 25 & 125 \\
        & & & \\
        1 & 8 & 64 & 512 \\
      \end{array}
    \right).
$$

For this matrix we have

$$
B = \mathcal{BD}(A) = \left(
      \begin{array}{cccc}
        1 & 2 & 2 & 2 \\
        & & & \\
        1 & 1 & 3 & 3 \\
        & & & \\
        1 & 2 & 6 & 5 \\
        & & & \\
        1 & 3/2 & 5/2 & 90 \\
      \end{array}
    \right),
$$

which means

$$
A = {\bar F}_3{\bar F}_2{\bar F}_1 D {\bar G}_1{\bar G}_2{\bar G}_3 =
$$

$$
= \left(
      \begin{array}{cccc}
        1 & 0 & 0 & 0 \\
        & & & \\
        0 & 1 & 0 & 0 \\
        & & & \\
        0 & 0 & 1 & 0 \\
        & & & \\
        0 & 0 & 1 & 1 \\
      \end{array}
    \right)
    \left(
      \begin{array}{cccc}
        1 & 0 & 0 & 0 \\
        & & & \\
        0 & 1 & 0 & 0 \\
        & & & \\
        0 & 1 & 1 & 0 \\
        & & & \\
        0 & 0 & \frac{3}{2} & 1 \\
      \end{array}
    \right)
    \left(
      \begin{array}{cccc}
        1 & 0 & 0 & 0 \\
        & & & \\
        1 & 1 & 0 & 0 \\
        & & & \\
        0 & 2 & 1 & 0 \\
        & & & \\
        0 & 0 & \frac{5}{2} & 1 \\
      \end{array}
    \right)
    \left(
      \begin{array}{cccc}
        1 & 0 & 0 & 0 \\
        & & & \\
        0 & 1 & 0 & 0 \\
        & & & \\
        0 & 0 & 6 & 0 \\
        & & & \\
        0 & 0 & 0 & 90 \\
      \end{array}
    \right)
$$
$$
    \left(
      \begin{array}{cccc}
        1 & 2 & 0 & 0 \\
        & & & \\
        0 & 1 & 3 & 0 \\
        & & & \\
        0 & 0 & 1 & 5 \\
        & & & \\
        0 & 0 & 0 & 1 \\
      \end{array}
    \right)
    \left(
      \begin{array}{cccc}
        1 & 0 & 0 & 0 \\
        & & & \\
        0 & 1 & 2 & 0 \\
        & & & \\
        0 & 0 & 1 & 3 \\
        & & & \\
        0 & 0 & 0 & 1 \\
      \end{array}
    \right)
    \left(
      \begin{array}{cccc}
        1 & 0 & 0 & 0 \\
        & & & \\
        0 & 1 & 0 & 0 \\
        & & & \\
        0 & 0 & 1 & 2 \\
        & & & \\
        0 & 0 & 0 & 1 \\
      \end{array}
    \right),
$$
and

$$
A^{-1} = G_1G_2G_3D^{-1}F_3F_2F_1.
$$

$$
A^{-1}= \left(
      \begin{array}{cccc}
        1 & -2 & 0 & 0 \\
        & & & \\
        0 & 1 & -2 & 0 \\
        & & & \\
        0 & 0 & 1 & -2 \\
        & & & \\
        0 & 0 & 0 & 1 \\
      \end{array}
    \right)
    \left(
      \begin{array}{cccc}
        1 & 0 & 0 & 0 \\
        & & & \\
        0 & 1 & -3 & 0 \\
        & & & \\
        0 & 0 & 1 & -3 \\
        & & & \\
        0 & 0 & 0 & 1 \\
      \end{array}
    \right)
    \left(
      \begin{array}{cccc}
        1 & 0 & 0 & 0 \\
        & & & \\
        0 & 1 & 0 & 0 \\
        & & & \\
        0 & 0 & 1 & -5 \\
        & & & \\
        0 & 0 & 0 & 1 \\
      \end{array}
    \right)
$$
$$
    \left(
      \begin{array}{cccc}
        1 & 0 & 0 & 0 \\
        & & & \\
        0 & 1 & 0 & 0 \\
        & & & \\
        0 & 0 & 1/6 & 0 \\
        & & & \\
        0 & 0 & 0 & 1/90 \\
      \end{array}
    \right)
    \left(
      \begin{array}{cccc}
        1 & 0 & 0 & 0 \\
        & & & \\
        0 & 1 & 0 & 0 \\
        & & & \\
        0 & 0 & 1 & 0 \\
        & & & \\
        0 & 0 & \frac{-5}{2} & 1 \\
      \end{array}
    \right)
    \left(
      \begin{array}{cccc}
        1 & 0 & 0 & 0 \\
        & & & \\
        0 & 1 & 0 & 0 \\
        & & & \\
        0 & -2 & 1 & 0 \\
        & & & \\
        0 & 0 & \frac{-3}{2} & 1 \\
      \end{array}
    \right)
    \left(
      \begin{array}{cccc}
        1 & 0 & 0 & 0 \\
        & & & \\
        -1 & 1 & 0 & 0 \\
        & & & \\
        0 & -1 & 1 & 0 \\
        & & & \\
        0 & 0 & -1 & 1 \\
      \end{array}
    \right).
$$

On the other hand, the factorization of $A^{-1}$ corresponding to the matrix interpretation of the Bj\"orck-Pereyra algorithm is the following (see \cite{GolV}):

$$
A^{-1} =  G_1G_2G_3D_2^{-1}L_2D_1^{-1}L_1D_0^{-1}L_0 =
$$

$$
= G_1G_2G_3
\left(
      \begin{array}{cccc}
        1 & 0 & 0 & 0 \\
        & & & \\
        0 & 1 & 0 & 0 \\
        & & & \\
        0 & 0 & 1 & 0 \\
        & & & \\
        0 & 0 & 0 & \frac{1}{6} \\
      \end{array}
    \right)
    \left(
      \begin{array}{cccc}
        1 & 0 & 0 & 0 \\
        & & & \\
        0 & 1 & 0 & 0 \\
        & & & \\
        0 & 0 & 1 & 0 \\
        & & & \\
        0 & 0 & -1 & 1 \\
      \end{array}
    \right)
    \left(
      \begin{array}{cccc}
        1 & 0 & 0 & 0 \\
        & & & \\
        0 & 1 & 0 & 0 \\
        & & & \\
        0 & 0 & \frac{1}{3} & 0 \\
        & & & \\
        0 & 0 & 0 & \frac{1}{5} \\
      \end{array}
    \right)
$$
$$
    \left(
      \begin{array}{cccc}
        1 & 0 & 0 & 0 \\
        & & & \\
        0 & 1 & 0 & 0 \\
        & & & \\
        0 & -1 & 1 & 0 \\
        & & & \\
        0 & 0 & -1 & 1 \\
      \end{array}
    \right)
    \left(
      \begin{array}{cccc}
        1 & 0 & 0 & 0 \\
        & & & \\
        0 & 1 & 0 & 0 \\
        & & & \\
        0 & 0 & \frac{1}{2} & 0 \\
        & & & \\
        0 & 0 & 0 & \frac{1}{3} \\
      \end{array}
    \right)
    \left(
      \begin{array}{cccc}
        1 & 0 & 0 & 0 \\
        & & & \\
        -1 & 1 & 0 & 0 \\
        & & & \\
        0 & -1 & 1 & 0 \\
        & & & \\
        0 & 0 & -1 & 1 \\
      \end{array}
    \right)
$$

\medskip

Let us observe that Stage II of the algorithms (corresponding to the product $G_1G_2G_3$) is the same for both factorizations. This matrix $G_1G_2G_3$ is the matrix of change of basis from the Newton basis to the monomial basis.

On the contrary, Stage I of the algorithm (corresponding to the factorizations $D^{-1}F_3F_2F_1 = D_2^{-1}L_2D_1^{-1}L_1D_0^{-1}L_0$) is different in the the Bj\"orck-Pereyra algorithm and in the algorithm {\tt TNSolve} which starts from $B = \mathcal{BD}(A)$. This first stage corresponds to the computation of the {\sl divided differences}, i.e. the coefficients of the interpolating polynomial in the Newton basis.

\medskip

We finish this section with some numerical experiments of linear system solving with Vandermonde and Cauchy matrices. We only include simple experiments with matrices of small size, and using available algorithms, so that the interested reader can easily reproduce them. All the algorithms of Plamen Koev, put together in the package {\tt TNTool},  are available in \cite{KOEV}.

We begin with the Vandermonde case, with a Vandermonde matrix $A$ of order $n = 7$ corresponding to the nodes $1,2,3,4,5,6,7$ (these nodes are the second column of the matrix), whose condition number is $\kappa_2(A) = 2.4e+7$. We will solve the linear system $A x = f$, with
$$
f = [1/21,-1/21,1/23,-1/23,1/29,-1/29,1/31]^T.
$$

We will compare the approximate solutions with the exact solution $x_e$ computed in exact rational arithmetic by using {\sl Maple}, and the relative errors will be computed in MATLAB as ${\tt norm(x - xe,2)/norm(xe,2)}$, i.e. we are computing the relative error in the Euclidean norm.

First, we compute the solution by means of the MATLAB function ${\tt A \backslash f}$, and the corresponding relative error is $1.5e-13$. Next, the system is solved by the Bj\"orck-Pereyra algorithm by using the algorithm {\tt VTsolve}, which can be obtained from the m-files of Chapter $4$ of the book of Golub and Van Loan \cite{GolV}, available in its web page: {\tt www.cs.cornell.edu/cv/GVL4/golubandvanloan.htm}. In this case the relative error is $2.6e-16$, which confirms the high relative accuracy to be expected of this algorithm when $f$ has an alternating sign pattern.

Finally, we compute the solution of the linear system by using the algorithm {\tt TNSolve} of P. Koev, by previously computing $B = \mathcal{BD}(A)$ by means of the algorithm {\tt TNVandBD} of P. Koev. In this case the relative error is $2.6e-16$, which again confirms the high relative accuracy to be expected of this algorithm when $f$ has an alternating sign pattern.

\medskip

In the second example the coefficient matrix is a Hilbert matrix $A$ of order $n = 7$, constructed in MATLAB by means of the instruction $\tt A = hilb(7)$, whose condition number is $\kappa_2(A) = 4.7e+8$. We will solve the linear system $A x = f$, with
$$
f = [1/21,-1/21,1/23,-1/23,1/29,-1/29,1/31]^T.
$$

As before,we will compare the approximate solutions with the exact solution $x_e$ computed in exact rational arithmetic by using {\sl Maple}, and the relative errors will be computed in MATLAB as ${\tt norm(x - xe,2)/norm(xe,2)}$.

First, we compute the solution by means of the MATLAB function ${\tt A \backslash f}$, and the corresponding relative error is $5.0e-9$. Next, the system is solved by the {\tt BKO} algorithm of Boros, Kailath and Olshevsky, which can be taken from page $277$ of their paper \cite{BKO}. In this case the relative error is $1.3e-16$, which confirms the high relative accuracy to be expected of this algorithm when $f$ has an alternating sign pattern.

Finally, we compute the solution of the linear system by using the algorithm {\tt TNSolve} of P. Koev, by previously computing $B = \mathcal{BD}(A)$ by means of the algorithm {\tt TNCauchyBD} of P. Koev, with the instruction {\tt B = TNCauchyBD(x,y)}, by defining {\tt x = [0,1,2,3,4,5,6]} and {\tt y = [1,2,3,4,5,6,7]}. In this case the relative error is $1.4e-16$, which again confirms the high relative accuracy to be expected of this algorithm when $f$ has an alternating sign pattern.

\section{Eigenvalue and singular value problems}

In Section $59.3$ of the already commented chapter \cite{Drmac}, Z. Drma{\v c} refers to the work of Koev \cite{KOEV05,KOEV07} to recall that if the totally nonnegative matrix $A$ is given implicitly by its bidiagonal decomposition the all its singular values (and eigenvalues, too) can be computed to high relative accuracy. The author also recalls that an accurate bidiagonal representation is possible provided certain minors can be accurately computed.

These results are not included in the recent books of Bj\"orck \cite{Bjorck} and Golub-Van Loan \cite{GolV}, but is important to indicate that {\it the same bidiagonal factorization} (with $n^2$ parameters stored in the matrix $B = \mathcal{BD}(A)$) is the starting point for accurately computing the eigenvalues and the singular values of $A$. By using the algorithms $\tt TNEigenValues$ and $\tt TNSingularValues$ of the package $\tt TNTool$ \cite{KOEV}, one can compute the eigenvalues and the singular values (respectively) of $A$ to high relative accuracy.

The application of the bidiagonal decomposition for computing eigenvalues and singular values of totally positive matrices developed by P. Koev \cite{KOEV05,KOEV07} has an important precedent in the work of Demmel and coworkers \cite{DGE}. In section $9$ of \cite{DGE}, devoted to the class of totally positive matrices, the authors indicate that ``{\it achieving high relative accuracy requires not just total positivity but an appropriate parameterization that permits minors to be evaluated to high relative accuracy}''.

Some years later, in section $2$ of \cite{KOEV05}, the author acknowledges these contributions of \cite{DGE}, and in section $3$ he adds a fundamental tool: the results on Neville elimination and total positivity introduced by Gasca and Pe\~na in \cite{GP92,GP94,GP96}. Koev shows that being able to compute all $n^2$ initial minors is a necessary and sufficient condition for accurately computing the bidiagonal decomposition ($\mathcal{BD}(A)$) of a totally positive matrix.

In section $7$ of \cite{KOEV05} we find the main idea for the construction of accurate algorithms: ``{\it In other words, $\mathcal{BD}(A)$ determines the eigenvalues and the singular values of $A$ accurately''}. In addition, it is indicated that the final step of the algorithms $\tt TNEigenValues$ and $\tt TNSingularValues$ is the computation of the singular values of a bidiagonal matrix by using the LAPACK routine $\tt DLASQ1$ (\cite{DKah,FP}) (which must be compiled to be used in the MATLAB algorithms), and it is known to introduce only a small additional relative error (\cite{FP}).

\bigskip

As we read in the introduction of \cite{KOEV05}, {\sl``when traditional algorithms are used to compute the eigenvalues or the singular values of an ill-conditioned TN matrix, only the largest eigenvalues and the largest singular values are computed with guaranteed relative accuracy. The tiny eigenvalues and singular values may be computed with no relative accuracy at all, even though they may be the only quantities of practical interest'' }.
We will illustrate with two simple examples the good behaviour of the algorithms {\tt TNEigenValues} and {\tt TNSingularValues} when applied to ill-conditioned totally positive matrices.

\medskip

The first matrix we consider is a Hilbert matrix of order $10$, constructed in MATLAB by means of the instruction $\tt A = hilb(10)$. Let us recall that a Hilbert matrix is a special case of Cauchy matrix with generic entries $c_{ij} = 1/(x_i - y_j)$. The condition number of this matrix is $\kappa_2(A) = 1.6e+13$.

We can compute $B = \mathcal{BD}(A)$ by means of the algorithm {\tt TNCauchyBD} of P. Koev, with the instruction {\tt B = TNCauchyBD(x,y)}, by defining
$$
{\tt x = [0,1,2,3,4,5,6,7,8,9]}, {\tt y = [1,2,3,4,5,6,7,8,9,10]}
$$
(P. Koev is using as entries of the Cauchy matrix $c_{ij} = 1/(x_i + y_j)$). Then the eigenvalues are computed by means of the instruction {\tt TNEigenvalues(B)}.

Then these eigenvalues are also computed by means of the standard MATLAB instruction {\tt eig(A)}. We compare the approximate value of the smallest eigenvalue with the ``exact'' value obtained by using {\sl Maple} with extended precision, and we obtain for the value computed by means of {\tt eig(A)} a relative error of $2.3e-4$, while for the algorithm {\tt TNEigenValues} the relative error is $3.4e-16$.

\medskip

In the  second example we compute the singular values of the Pascal matrix of order $10$, constructed in MATLAB by means of the instruction $\tt A = pascal(10)$. Now the condition number is $\kappa_2(A) = 4.1e+9$.

It is interesting to observe that for Pascal matrices we have an exact $B = \mathcal{BD}(A)$: the matrix with all the entries equal to $1$, which is constructed in MATLAB as {\tt B = ones(n,n)}. Then the singular values are computed by means of the instruction {\tt TNSingularvalues(B)}.

These eigenvalues are also computed by means of the standard MATLAB instruction {\tt svd(A)}. We compare the approximate value of the smallest singular value with the ``exact'' value obtained by using {\sl Maple} with extended precision, and we obtain for the value computed by means of {\tt eig(A)} a relative error of $2.9e-9$, while for the algorithm {\tt TNEigenValues} the relative error is $6.5e-16$.

\bigskip
We see the effect of the ill-conditioning of the Hilbert and Pascal matrices when using the standard MATLAB functions {\tt eig} and {\tt svd}, while the algorithms $\tt TNEigenValues$ and $\tt TNSingularValues$ (starting from an accurate $\mathcal{BD}(A)$) give high relative accuracy.

The availability of the algorithms of Koev in \cite{KOEV} has encouraged the search for new algorithms for the bidiagonal decomposition of various classes of totally positive structured matrices, including matrices from new application fields such as \cite{DPP} (see \cite{PenaS} for a recent account).

\section{The rectangular case}

Linear system solving and eigenvalue computation correspond to the case in which $A$ is a square matrix but, as seen in \cite{KOEV07}, this bidiagonal decomposition also exists for rectangular matrices. For instance, an algorithm for computing it for the case of Bernstein-Vandermonde matrices is presented in \cite{MM13} (see algorithm {\tt TNBDBVR} in \cite{KOEV}).

Let us recall that in the square case the matrices must be nonsingular. For the extension to the rectangular case we must take into account the following comment in the Introduction of \cite{KOEV07}:

``The existence and uniqueness of the bidiagonal decomposition is critical to the design of our algorithms. Therefore we restrict the class of totally nonnegative matrices under consideration to only those that are {\sl leading contiguous submatrices of square nonsingular totally nonnegative matrices}''.

These rectangular matrices arise in a natural way when solving least squares problems. In this context, and returning to our guide book \cite{Handbook}, we find this subject studied in Chapter $52$, entitled {\sl Least Squares Solution of Linear Systems} \cite{Hansen}. In sections $52.4$ and $52.5$, mainly following the classical book of Bj\"orck \cite{BJO}, the authors recall the use of the $QR$ factorization of $A$ (of size $m \times n$ with $m \ge n$), which leads to solving a (square) linear system $R x = Q^T b$. As read in \cite{Hansen}, the algorithm that is least sensitive to the influence of rounding errors is based on the $QR$ factorization of $A$, as first suggested by G. H. Golub in \cite{Gol65}.

For this purpose, P. Koev includes in \cite{KOEV07} a new gem for computing with totally positive matrices: the $QR$ factorization of a totally positive matrix $A$. Starting from B = $\mathcal{BD}(A)$ (the bidiagonal decomposition of $A$), the algorithm {\tt TNQR} computes $Q$ and $\mathcal{BD}(R)$ (the bidiagonal decomposition of the triangular factor $R$). This algorithm has been used in \cite{MM15} to solve least squares problems by using the Bernstein basis, including the accurate computation of the projection matrix (the {\sl hat matrix} of statistics). An extension to the bivariate setting has recently been carried out in \cite{MMV}.

Another important situation where the bidiagonal factorization for the rectangular case is used is the computation of singular values (one of its applications being the computation of the spectral condition number of a matrix), as seen in section $7$ of \cite{MM13}, and also the computation of the Moore-Penrose inverse, as carried out in \cite{MM18}.

\section{The tridiagonal case}

In this section we briefly consider the case in which the totally positive matrix $A$ is also symmetric and tridiagonal. An important example of symmetric and tridiagonal matrices are the Jacobi matrices associated with a family of orthogonal polynomials, and in this context the eigenvalues of the Jacobi matrices are the zeros of the corresponding orthogonal polynomials (\cite{Gautschi}, \cite{M2001}), which consequently can be computed with high relative accuracy starting from an accurate $\mathcal{BD}(A)$.

An example of this situation (where the Jacobi matrices are totally positive) has been analyzed in \cite{MM15a}.

An interesting aspect of this symmetric and tridiagonal case is the fact that now Neville elimination is the same as Gaussian elimination, and the bidiagonal decomposition of $A$ is precisely $A = LDL^T$, where $L$ is now not only lower triangular but also bidiagonal. Although the author does not pay explicit attention to total positivity, the algorithmic advantages of using the $LDL^T$ factorization have been shown by B. N. Parlett in \cite{Parlett}.

\section{Conclusions}

The starting point of our work has been the paper of Bj\"orck and Pereyra of $1970$ \cite{BP}, where they showed that the Newton-Horner algorithm for solving a Vandermonde system associated to polynomial interpolation could be interpreted in connection with a factorization of the inverse of the Vandermonde matrix as a product of bidiagonal matrices.

Our main goal has been to show how this approach to solve linear systems has been extended to solve several other linear algebra problems, such as eigenvalue and singular value computation, or least squares problems. The main tool is the use of the bidiagonal decomposition of the corresponding matrix (not only of its inverse, as in the Bj\"orck-Pereyra algorithm), stored as $B = \mathcal{BD}(A)$.

This idea, first presented in the work of Koev \cite{KOEV05} (see also \cite{DDHK}), along with the availability of his algorithms, has encouraged the search for new algorithms for the bidiagonal decomposition of various classes of totally positive structured matrices.

It must be remarked that although a superficial reading (or a careless writing) of our references could suggest that Neville elimination is the method to construct those algorithms, the fact is that Neville elimination is the fundamental {\sl theoretical tool} but not the algorithmic tool (see section $5$ in \cite{MM18}). For achieving high relative accuracy the algorithms must be adapted to the specific structure of each class of matrices.

\begin{acknowledgements}
This research has been partially supported by Spanish Research Grant MTM2015-65433-P (MINECO/FEDER) from the Spanish Ministerio de Econom{\'\i}a y Competitividad. J. J. Mart{\'\i}nez is member of the Research Group {\sc asynacs} (Ref. {\sc ccee2011/r34)} of Universidad de Alcal\'a.

\end{acknowledgements}



\begin{thebibliography}{1}

\bibitem{Bella}
Bella, T., Eidelman, Y., Gohberg, I., Koltracht, I., Olshevsky, V.:
\newblock A Bj\"orck-Pereyra-type algorithm for Szeg\"o-Vandermonde matrices based on properties of unitary Hessenberg matrices.
\newblock {\em Linear Algebra and Its Applications}, 420:634--647, 2007.

\bibitem{BJO}
Bj\"orck, A.:
\newblock {\em Numerical Methods for Least Squares Problems}.
\newblock SIAM, Philadelphia, 1996.

\bibitem{Bjorck}
Bj\"orck, A.: \textit{Numerical Methods in Matrix Computations}. Texts in Applied Mathematics, Volume $59$. Springer International Publishing, Switzerland, 2015.

\bibitem{BP}
Bj\"orck, A., Pereyra, V.: Solution of Vandermonde systems of equations. \textit{Math. Comp.} 24:893-903, 1970.

\bibitem{BKO}
Boros, T., Kailath, T., Olshevsky, V.: A fast parallel Bj\"orck-Pereyra-type algorithm for solving Cauchy linear equations. \textit{Linear Algebra Appl.} 302/303:265-293, 1999.

\bibitem{DPP}
Delgado, J., Pe\~na, G, Pe\~na, J. M.:
\newblock Accurate and fast computations with positive extended Schoenmakers-Coffey matrices.
\newblock {\em Numer. Linear Algebra Appl.}, 23:1023--1031, 2016.


\bibitem{DGE}
Demmel, J., Gu, M., Eisenstat, S., Slapni{\v c}ar, I., Veseli\'c, K., Drma{\v c}, Z.: Computing the singular value decomposition with high relative accuracy. \textit{Linear Algebra and Its Applications} 299:21-80, 1999.

\bibitem{DDHK}
Demmel, J., Dumitriu, I., Holtz, O., Koev, P.:
\newblock Accurate and efficient expression evaluation and linear algebra.
\newblock {\em Acta Numerica}, 17:87--145, 2008.




\bibitem{DKah}
Demmel, J., Kahan, W.: Accurate singular values of bidiagonal matrices. \textit{SIAM J. Sci. Stat. Comput} 11(5):873--912, 1990.


\bibitem{DK}
Demmel, J., Koev, P.: The accurate and efficient solution of a totally positive generalized Vandermonde linear system. {\em SIAM Journal on Matrix Analysis and Applications}, 27(1):142--152, 2005.


\bibitem{Drmac}
Drma{\v c}, Z.: Computing Eigenvalues and Singular Values to High Relative Accuracy. Chapter $59$ in: L. Hogben (ed.), \textit{Handbook of Linear Algebra (Second Edition)}.  CRC Press, Boca Raton (FL), 2014.

\bibitem{FJ}
Fallat, S. M., Johnson, C. R.:
\newblock {\em Totally Nonnegative Matrices}.
\newblock  Princeton University Press, Princeton and Oxford, 2011.


\bibitem{FP}
Fernando, K., Parlett, B.: Accurate singular values and differential qd algorithms. \textit{Numerische Mathematik} 67:191-229, 1994.

\bibitem{GMi}
Gasca, M., Micchelli, C. A. (eds.):
\newblock {\em Total Positivity and Its Applications (Jaca, 1994)}. Kluwer Academic Publishers, Dordrecht, 1996.



\bibitem{GP92}
Gasca, M., Pe\~na, J. M.:
\newblock Total positivity and Neville elimination.
\newblock {\em Linear Algebra and Its Applications}, 165:25--44, 1992.


\bibitem{GP94}
Gasca, M., Pe\~na, J. M.:
\newblock A matricial description of Neville elimination with applications to total positivity.
\newblock {\em Linear Algebra and Its Applications}, 202:33--45, 1994.

\bibitem{GP96}
Gasca, M., Pe\~na, J. M.:
\newblock On Factorizations of Totally Positive Matrices.
\newblock In: M. Gasca and C. A. Michelli (Eds.), Total Positivity and Its Applications, Kluwer Academic Publishers, Dordrecht, 1996, pp. 109--130.

\bibitem{Gautschi}
Gautschi, W.:
\newblock {\em Orthogonal Polynomials. Computation and Approximation}.
\newblock Oxford University Press, Oxford, 2004.


\bibitem{Gol65}
Golub, G. H.: Numerical methods for solving least squares problems. \textit{Numerische Mathematik} 7:206--216, 1965.



\bibitem{GolV}
Golub, G. H., Van Loan, C. F.: \textit{Matrix Computations, 4th edition}. Johns Hopkins University Press, Baltimore, 2013.



\bibitem{Hansen}
Hansen, P. C., Nielsen, H. B.: Least Squares Solution of Linear Systems. Chapter $52$ in: L. Hogben (ed.), \textit{Handbook of Linear Algebra (Second Edition)}.  CRC Press, Boca Raton (FL), 2014.



\bibitem{Higham}
Higham, N. J.: Error analysis of the Bj\"orck-Pereyra algorithms for solving Vandermonde systems. \textit{Numer. Math.} 50:613-632, 1987.


\bibitem{HIGSIAM}
Higham, N. J.:
\newblock Stability analysis of algorithms for solving confluent Vandermonde-like systems.
\newblock {\em SIAM Journal on Matrix Analysis and Applications}, 11:23--41, 1990.



\bibitem{HIG}
Higham, N. J.:
\newblock  {\em Accuracy and Stability of Numerical Algorithms, second edition}.
\newblock SIAM, Philadelphia, 2002.

\bibitem{Handbook}
Hogben, L. (editor):
\newblock {\em Handbok of Linear Algebra (Second Edition)}. CRC Press, Boca Raton (FL), 2014.


\bibitem{KOEV}
Koev, P.:
\newblock {\tt http://www.math.sjsu.edu/$\sim$koev/}

\bibitem{KOEV05}
Koev, P.:
\newblock Accurate eigenvalues and SVDs of totally nonnegative matrices.
\newblock {\em SIAM Journal on Matrix Analysis and Applications}, 27:1--23, 2005.

\bibitem{KOEV07}
Koev, P.:
\newblock Accurate computations with totally nonnegative matrices.
\newblock {\em SIAM Journal on Matrix Analysis and Applications}, 29:731--751, 2007.



\bibitem{MM13}
Marco, A., Mart{\'\i}nez, J.-J.:
\newblock Accurate computations with totally positive Bernstein-Vandermonde matrices.
\newblock {\em Electron. J. Linear Algebra}, 26:357--380, 2013.


\bibitem{MM15a}
Marco, A., Mart{\'\i}nez, J.-J.:
\newblock A total positivity property of the Marchenko-Pastur law.
\newblock {\em Electron. J. Linear Algebra}, 30:106--117, 2015.


\bibitem{MM15}
Marco, A., Mart{\'\i}nez, J.-J.:
\newblock Ajuste polin\'omico por m{\'\i}nimos cuadrados usando la base de Bernstein.
\newblock {\em La Gaceta de la RSME}, 18:135--153, 2015.

\bibitem{MM18}
Marco, A., Mart{\'\i}nez, J.-J.:
\newblock Accurate computation of the Moore-Penrose of strictly totally positive matrices.
\newblock {\em Journal of Computational and Applied Mathematics} (in press), 2018.


\bibitem{MMP}
Marco, A., Mart{\'\i}nez, J.-J., Pe\~na, J. M.:
\newblock Accurate bidiagonal decomposition of totally positive Cauchy-Vandermonde matrices and applications.
\newblock {\em Linear Algebra and Its Applications}, 517:63--84, 2017.


\bibitem{MMV}
Marco, A., Mart{\'\i}nez, J.-J., Via\~na, R.:
\newblock Least squares problems involving generalized Kronecker products and application to bivariate polynomial regression.
\newblock {\em Numerical Algorithms} (in press), 2018.


\bibitem{M2001}
Mart{\'\i}nez, J.-J.:
\newblock Polinomios ortogonales, cuadratura gaussiana y problemas de valores propios.
\newblock In {\em Margarita Mathematica en memoria de Jos\'e Javier (Chicho) Guadalupe Hern\'andez}, (L. Espa\~nol and J. L. Varona,eds.), pp. 595--606, Servicio de Publicaciones, Universidad de La Rioja, Logro\~no, Spain, 2001.


\bibitem{MP98}
Mart{\'\i}nez, J.-J., Pe\~na, J. M.:
\newblock Factorizations of Cauchy-Vandermonde matrices.
\newblock {\em Linear Algebra and Its Applications}, 284:229--237, 1998.

\bibitem{MG87}
M\"uhlbach, G., Gasca, M.:
\newblock A test for strict total positivity via Neville elimination.
\newblock In {\em Current trends in Matrix Theory (Auburn, Ala., 1986)}, (F. Uhlig and R. Grone,eds.), pp. 225--232, North Holland, Amsterdam, 1987.




\bibitem{Parlett}
Parlett, B. N.:
\newblock For tridiagonals $T$ replace $T$ with $LDL^T$.
\newblock {\em J. Comput. Appl. Math.}, 123:117--130, 2000.


\bibitem{PenaS}
Pe\~na, J. M.: Accurate Computations and Applications of Some Classes of Matrices. In:
\newblock Accurate Computations and Applications of Some Classes of Matrices.
\newblock In: Mateos, M., Alonso, P. (eds.) {\em Computational Mathematics, Numerical Analysis and Applications}. SEMA SIMAI Springer Series, vol. 13. Springer, Cham, 2017.

\bibitem{Pinkus}
Pinkus, A.:
\newblock {\em Totally Positive Matrices}.
\newblock Cambrige Tracts in Mathematics, Num. 181, Cambrigde University Press, 2010.



\bibitem{Stewart}
Stewart, M.: Fast Algorithms for Structured Matrix Computations. Chapter $62$ in: L. Hogben (ed.), \textit{Handbook of Linear Algebra (Second Edition)}.  CRC Press, Boca Raton (FL), 2014.





\end{thebibliography}
\end{document}